\documentclass[twoside]{article}
\usepackage{amssymb,a4}
\usepackage{amsthm}
\def\cen{\centerline}
\usepackage{indentfirst}
\def\pt{\partial}

\def\z{\mathbb{Z}}
\def\c{\mathbb{C}}

\def\ad{\hbox{ad}}

\def\a{\alpha}
\def\b{\beta}
\def\sg{\sigma}

\def\mod{\hbox{mod}}
\def\gl{\frak{gl}}
\newfont{\df}{eufm10}
\def\vep{\varepsilon}

\newcommand{\fr}[2]{\frac{\displaystyle #1}{\displaystyle #2}}

\def\sg{\sigma}

\def\der{\hbox{Der}}
\def\stl{\frak{stl}}
\def\fsl{\frak{sl}}

\def\r{\gamma}

\def\ot{\otimes}

\def\de{\delta}

\def\im{{\hbox{\rm Im}\,}}

\def\hom{{\hbox{\rm Hom}}}
\def\aut{{\hbox{\rm Aut}\,}}

\def\ad{{\hbox{\rm ad}\,}}

\def\ot{\otimes}
\def\wg{\tilde{\frak g}}
\def\mg{{\bf \frak g}}
\def\dmg{\dot{\mg}}

\hoffset
\voffset
\title{{ Leibniz Superalgebras and Central Extensions}
\footnotetext{* Corresponding author. Email: ldecnu2001@yahoo.com,
liud@czu.cn}}
\author{{\bf Dong Liu}* \\ Department of Mathematics, Shanghai Jiaotong University\\ Shanghai, 200240, P.R. China
\\Department of Mathematics, Changzhou Institute of Technology\\ Changzhou, 213002, P.R. China\\ \\{\bf Naihong Hu}\\Department of Mathematics, East China Normal University\\ Shanghai, 200062, P.R. China }
\date{ }
\begin{document}
\maketitle
\def\abstractname{ABSTRACT}
\begin{abstract}
Dialgebras are generalizations of associative algebras which give
rise to Leibniz algebras instead of Lie algebras. In this paper we
study super dialgebras and Leibniz superalgebras, which are
$\z_2$-graded dialgebras and Leibniz algebras. We also study
universal central extesions of Leibniz superalgebras and obtain
some results as in the case of Lie superalgebras. We determine
universal central extensions of basic classical Lie superalgebras
in the category of Leibniz superalgebras. They play a key role in
studying Leibniz superalgebras graded by finite root systems.

\smallskip
\noindent
{\it Key Words}: Leibniz superalgebra; super dialgebra;
universal central extension.
\end{abstract}

\newtheorem{theo}{Theorem}[section]
\newtheorem{defi}[theo]{Definition}
\newtheorem{lemm}[theo]{Lemma}
\newtheorem{prop}[theo]{Proposition}
\newtheorem{coro}[theo]{Corollary}

\section{INTRODUCTION}

Central extensions play an important role in the theory of Lie
algebras, and it is therefore not surprising that there are many
results on central extensions of various classes of Lie algebras.
Recently several authors have considered central extensions of
Leibniz algebras (\cite{Lo2}, \cite{Gao1}, \cite{Gao2}, \cite{L},
\cite{LH1}, etc. ) and Lie superalgebras (\cite{IK}, \cite{MP},
etc.).

Universal central extensions of simple Lie algebras over
commutative rings were studied in \cite{Gar} and \cite{Kas}. To be
precise, let $\dmg$ be a simple finite dimensional Lie algebra
over the field of complex numbers $\c$. It is known that the
kernel of the universal central extension of the Lie algebra
$\dmg\ot A$, where $A=\c[t_1^\pm,\cdots,t_\nu^\pm]\,(\nu\ge1)$, is
$\Omega_A^1/dA$, the $A$-module of K$\ddot{\hbox{a}}$hler
differentials over $\c$ (\cite{Gar}, \cite{Kas}, \cite{KL},
\cite{MRY}). Some authors (\cite{Gao1}, \cite{Gao2}, \cite{LL},
etc.) studied the universal central extension of $\dmg\ot A$ in
the category of Leibniz algebras.  The kernel of this  universal
central extension of $\dmg\ot A$ is $\Omega_A^1$. The universal
central extensions of basic classical Lie superalgebras over an
associative algebra have also been studied (\cite{IK}, \cite{MP},
etc.). Let $\mg$ be a basic classical simple Lie superalgebra
(\cite{Kac}) over $\c$, then the universal central extension
$\tilde U$ of $\mg\ot A$ is given by \cite{IK} as follows:
$$\mg(A):=\bar\mg\ot A\oplus\Omega_A^1/dA, $$
where $$\bar\mg=\cases{\mg, \quad \hbox{if}\ \mg\ \hbox{is not of type}\ A(n, n)(\forall n)\cr
                \fsl(n+1, n+1), \quad \hbox{if}\ \mg\ \hbox{is of type}\ A(n, n)(n>1)\cr
                \frak d, \quad \hbox{if}\ \mg\ \hbox{is of type}\ A(1, 1)          }\leqno(1.1)$$
and $\frak d$ is defined as in \cite{IK}.

The concepts of Leibniz superalgebra and its cohomology were first
introduced by Dzhumadil'daev in \cite{D}. Moreover pre-simplicial
structure on the cochain complex of colour Leibniz algebras is
constructed in \cite{D}. The aim of this paper is to study super
dialgebras and Leibniz superalgebras. In Section 2 we recall some
notions of dialgebras and Leibniz algebras and their (co)homology
and Lie superalgebras. In Section 3 and Section 4 we define super
dialgebras and Leibniz superalgebras and their universal central
extensions. In Section 5 we mainly study the universal central
extensions of $\mg\ot D$ for a basic classical Lie superalgebra
$\mg$ in the category of Leibniz superalgebras over a unital
dialgebra $D$. Especially we obtain that the universal central
extension $\wg$ of $\mg\ot A$ is
$$\wg:=\bar\mg\ot A\oplus\Omega_A^1.$$

Throughout this paper $\c, \z, K$ denote the field of complex
numbers, the ring of integers and a field respectively, $|x|$
denotes the degree of a homogeneous element $x$ in a $\z_2$-graded
vector space, $\mg$ denotes a basic classical simple Lie
superalgebra.

\section{BASICS}

We recall some notions of dialgebras and Leibniz algebras and
their (co)homology as defined in \cite{Lo1}---\cite{Lo3} and some
notions of Lie superalgebras introduced in \cite{Kac} and
\cite{IK}.

\subsection{Dialgebras and Leibniz algebras}

\begin{defi} \cite{Lo2}  A (associative)
dialgebra $D$ over $K$ is a $K$-vector space equipped with two
operations $\dashv, \vdash:D\ot D\to D$, called left and right
products, satisfying the following five axioms:
$$\cases{a\dashv(b\dashv c)=(a\dashv b)\dashv c=a\dashv(b\vdash c),\cr
        (a\vdash b)\dashv c=a\vdash(b\dashv c),\cr
         (a\vdash b)\vdash c=a\vdash (b\vdash c)=(a\dashv b)\vdash c.  }\eqno(D)$$

\end{defi}

A dialgebra is called unital if it is given a specified bar-unit:
an element $1\in D$ which is a unit for the left and right
products only on the bar-side, that is $1\vdash a=a=a\dashv 1$,
for any $a\in D$. Denote by {\bf Dias} and {\bf Ass} the
categories of dialgebras and associative algebras over $K$
respectively. Then the category {\bf Ass} is a full subcategory of
{\bf Dias}.

Obviously a dialgebra is an associative algebra if and only if
$a\dashv b=a\vdash b=ab$.

The (co)homology of dialgebras can be seen in [F1] and [F2]. Now
we introduce the K$\ddot{\hbox{a}}$hler differential modules over
a dialgebra $D$ defined in \cite{L}.

For a commutative associative dialgebra $D$ over $K$ (i.e.
$a\dashv b=b\vdash a$), the module of differential $(\Omega_D^1,
d)$ of $D$ is defined in the following way. Let $\{a_i\}$ be any
basis for $D$ over $K$ and let $F$ be the free left $D$-module on
a basis $\{\tilde da_i\}$, where $\{\tilde da_i\}$ is some set
equipotent with $\{a_i\}$. We treat $F$ as a $2$-sided $D$-module
by setting $b\dashv(\tilde da)=(\tilde da)\vdash b$ and
$b\vdash(\tilde da)=(\tilde da)\dashv b$ for all $a, b\in D$. Let
$\tilde d: D\to F$ be the $K$-linear map $\sum c_i\dashv
a_i\mapsto\sum c_i\dashv\tilde da_i$ and let $N$ be the
$D$-submodule generated by the relations $\tilde d(a\star
b)-(\tilde da)\star b+a\star(\tilde da), a, b\in D$,
$\star=\dashv, \vdash$. Then $\Omega_D^1:=F/N$ and the canonical
quotient map $a\mapsto \tilde da+K$ is the differential map $d:
D\to \Omega_D^1$.

Up to evident isomorphism, $(\Omega_D^1, d)$ is characterized by
the property that for every $D$-module $M$ and every derivation
$f: D\to M$ there is a unique $D$-module map $g: \Omega_D^1\to M$
such that $f=g\circ d$. In this way $\der_{K}(D,
M)\cong\hom_D(\Omega_D^1, M)$.

A {\it Leibniz algebra} \cite{Lo2} $L$ is a vector space over a
field $K$ equipped with a $K$-bilinear map
$$[-,-]: L\times L\to L$$
satisfying the Leibniz identity
$$[[x, y], z]= [x, [y, z]]-[y, [x, z]], \quad \forall \;x, \,y, \,z\in L.\eqno(2.1)$$

Obviously, a Lie algebra is a Leibniz algebra. A Leibniz algebra
is a Lie algebra if and only if $[x, x]=0$ for all $x\in L$.

\subsection{Chevalley basis of basic classical simple Lie
superalgebras}

Let $\mg$ be a basic classical simple Lie superalgebra (see
\cite{Kac} in details), $A=(a_{ij})_{i, j\in I}$ for some fixed
index set $I$ be the Cartan matrix of $\mg$. Let
$D=\hbox{diag}(\vep_i)_{i\in I}$ and $B=(b_{ij})$ be diagonal and
symmetric matrices such that $A=DB$. Let $\Delta$ be its root
system. We fix a Dynkin diagram. Hence an odd simple root is
unique, and is denoted by $\a_{i_0}$.

For $\a=\sum_{i\in I}k_i\a_i\in \Delta$, we define
$$\vep_\a:=\cases{\fr{2}{(\a, \a)}, \quad (\a, \a)\ne 0\cr   \vep_{i_0}, \hfill (\a, \a)=0}$$
and $$H_{\a}:=\vep_\a\sum_{i\in I}k_i\vep_i^{-1}h_i,$$where $h_i$
is the $i$-th simple coroot(see \cite{IK} in detail).

By Lemma 2.23 in \cite{IK} we have
$$\a(H_\a)=\cases{2, \quad (\a, \a)\ne 0\cr 0, \quad (\a, \a)=0}.$$

We also set
$$\sg_\a:=\cases{-1, \quad \a\in -\Delta_{\bar1}^+\cr  1,\hfill \hbox{otherwise}},$$ and by definition we have
$$\sg_{-\a}=(-1)^{|\a|}\sg_\a.$$
For each $\a\in\Delta$, choose a root vector $X_\a\in \mg_\a$ such that
$$[x_\a, x_{-\a}]=\sg_\a H_\a,$$
then we have for these root vectors, we define $N_{\a, \b}\in\c(\a, \b\in\Delta)$ by
$$N_{\a, \b}:=\cases{\hbox{the coefficient of}\  X_{\a+\b}\ \hbox{in}\ [X_\a, X_\b], \quad \a+\b\in\Delta\cr   0,\quad \a+\b\notin \Delta }$$

The Chevalley basis of a basic classical Lie superalgebra was
determined in \cite{IK}.

\begin{theo}$\cite{IK}$

1. Let $\mg$ be a basic classical Lie superalgebra over $\c$. For each $\a\in\Delta$, we define $\sg_\a$ and $H_{\a}$ as above. Then there exist vectors $\{X_\a\in\mg_\a|\a\in\Delta \}$ such that

1) $[X_\a, X_{-\a}]=\sg_\a H_\a$,

2) $[X_\a, X_\b]=N_{\a, \b}X_{\a+\b}$, where the structure constants $\{N_{\a, \b}\}$ satisfy the following conditions:

i) If $\a\in \Delta_{\bar0}$ or $\b\in\Delta_{\bar0}$ (we
assume  $\a\in \Delta_{\bar0}$) and $\a+\b\in \Delta$, then
$$N_{\a, \b}^2=(p+1)^2.$$

ii) If $\a, \b\in \Delta_{\bar1}$ and $(\a, \a)\ne 0$ or $(\b, \b)\ne0$(we assume $(\a, \a)\ne 0$) and $\a+\b\in \Delta$, then
$$N_{\a, \b}^2=(p+1)^2.$$

iii) If $\a, \b\in \Delta_{\bar1}$ and $(\a, \a)=(\b, \b)=0$ and $\a+\b\in \Delta$, then
$$N_{\a, \b}^2=\b(H_\a)^2.$$
Where $p:=\max\{i\mid\b-i\a\in\Delta\}$.

2. Let $\{X_\a'|\a\in\Delta \}$ be another set of root vectors
satisfying the above conditions and  $\{N_{\a, \b}|\a, \b\in
\Delta\}'$ be corresponding structure constants. Then there exists
$\{u_\a|\a\in\Delta\}\subset\{\pm1\}$ such that $u_\a u_{-\a}=1$
and
$$N_{\a, \b}'=u_\a u_\b u_{\a+\b}^{-1}N_{\a, \b}, $$ for any $\a, \b\in \Delta$.
\end{theo}
In the sequel, we call the set $\{H_\a, X_\a|\a\in\Delta\}$ a {\it
Chevalley basis} of $\mg$ if it satisfies the condition in Theorem
2.2.

We also need to use the following Lemma, which is given in
\cite{IK}.

\begin{lemm}\cite{IK}
Let $\Delta$ be the root system of a basic classical Lie superalgebra $\mg$.
If $\a\in\Delta_{\bar1}^+$ such that $(\a, \a)=0$, then there exists $\r\in \Delta$ such that $\a(H_\r)=0.$

\end{lemm}

\section{SUPER DIALGEBRAS AND LEIBNIZ SUPERALGEBRAS}

\subsection{Super dialgebras.}

To study Leibniz superalgebras, we introduce super dialgebras.

\begin{defi}
 Let $K$ be a field. A {\it super dialgebra} over $K$ is a $\z_2$-graded $K$-vector space $D$ with two operations $\dashv,\, \vdash:D\ot D\to D$, called left and right products, satisfying the axiom (D)(section 2)
and $$D_{\sg}\dashv D_{\sg'},\  D_{\sg}\vdash D_{\sg'} \subset
D_{\sg+\sg'},\quad \forall\sg, \sg'\in\z_2.$$

The definitions of bar-unit and (homogeneous) morphisms of super
dialgebras are similar. Denote by {\bf SDias, SAss} the categories
of super dialgebras and associative superalgebras over $K$
respectively. Then the category {\bf SAss} is a full subcategory
of {\bf SDias}.
\end{defi}

\noindent{\bf Examples}. 1. Obviously an associative superalgebra
is a super dialgebra if $a\dashv b=a\vdash b=ab$.

2. {\it Super differential dialgebra.} Let $(A=A_{\bar0}\oplus
A_{\bar1}, d)$ be a differential associative super
algebra($|d|=\bar0$). So by hypothesis, $d(ab)=(da)b+adb$ and
$d^2=0$. Define left and right products on $A$ by the formulas
$$x\dashv y=xdy, \quad x\vdash y=(dx)y.$$
Then $A$ equipped with these two products is a super dialgebra.

3. {\it Tensor product.} If $D$ and $D'$ are two super dialgebras,
then the tensor product $D\ot D'$ is also a super dialgebra with
$$(a\ot a')\star(b\ot b')=(-1)^{|a'||b|}(a\star b)\ot (a'\star b')\eqno(3.1)$$ for $\star=\dashv, \vdash.$

For instance $M_{m+n}(D):=M_{m+n}(K)\ot D$ is a super dialgebra if
$D$ is a super dialgebra and $M_{m+n}(K)$ is the superalgebra of
all $(m+n)\times (m+n)$-matrices over $K$.

Similar to the result in \cite{Lo3}, we also have:

\begin{theo}
The free super dialgebra on a $\z_2$-graded vector space $V$ is
the dialgebra $Dias(V)=T(V)\ot V\ot T(V)$ equipped with the
induced $\z_2$-grading.\hfill $\rule[-.23ex]{1.0ex}{2.0ex}$
\end{theo}

A bimodule over a super dialgebra $D$, also called a
representation, is a $\z_2$-graded $K$-module $M$ equipped with
two linear maps
$$\vdash, \dashv: D\ot M\to M$$
$$\vdash, \dashv: M\ot D\to M$$satisfying the axioms (D) whenever they make sense and preserving $\z_2$-gradation.

For a super dialgebra $S$, let $S_{SAss}$ be the quotient of $S$
by the ideal generated by elements $x\dashv y-x\vdash y$ for all
$x, y\in S$. It is clear that $S_{SAss}$ is an associative
superalgebra. The canonical surmorphism $S\to S_{SAss}$ is
universal among the maps from $S$ to associative superalgebras. In
other words the associativization functor $(-)_{SAss}:$ {\bf
SDias}$\to${\bf SAss} is left adjoint to $inc:$ {\bf
SAss}$\to${\bf SDias}.

\subsection{Leibniz superalgebra}

\begin{defi} \cite{D}  A Leibniz superalgebra is a $\z_2$-graded vector space
$L=L_{\bar 0}\oplus L_{\bar 1}$ over a field $K$ equipped with a
$K$-bilinear map $[-,-]: L\times L\to L$ satisfying

$$[L_{\sg}, L_{\sg'}]\subset L_{\sg+\sg'},\quad \forall \sg, \sg'\in\z_2$$
and the Leibniz identity
$$[[a, b], c]= [a, [b, c]]-(-1)^{|a||b|}[b, [a, c]], \quad \forall \;a, \,b, \,c\in L.\eqno(3.2)$$

\end{defi}

Obviously, $L_{\bar0}$ is a Leibniz algebra. Moreover any Lie
superalgebra is a Leibniz superalgebra and any Leibniz algebra is
a trivial Leibniz superalgebra. A Leibniz superalgebra is a Lie
superalgebra if and only if
$$[\,a, b\,]+(-1)^{|a||b| }[b, a]=0, \quad \forall\; a, b\in L.$$

\noindent{\bf Examples.} 1. Let ${\mg}$ be a Lie superalgebra, $D$
be a unital commutative dialgebra, then $\mg\ot D$ with Leibniz
bracket $[x\ot a, y\ot b]=[x, y]\ot(a\vdash b)$ is a Leibniz
superalgebra. Let ${\mg}$ be a basic classical simple Lie
superalgebra, then $\tilde\mg=\bar\mg\ot D\oplus \Omega_D^1$ with
the bracket
$$[x\ot a, y\ot b]=[x, y]\ot(a\vdash b) +(x, y)b\dashv da, \quad \forall a, b\in D, x, y\in \bar\mg,\eqno(3.4)$$
$$[\Omega_D^1,\  \tilde\mg]=0\eqno(3.5)$$
is also a Leibniz superalgebra, where $\bar\mg$ defined in (1.1)
and $(-, -)$  is an even invariant bilinear form of $\bar\mg$. In
fact we shall prove that it is the universal central extension of
$\mg\ot D$ in Section 5.

2. Tensor product. Let $\mg$ be a Lie superalgebra,  then the bracket
$$[x\ot y, a\ot b]=[[x, y], a]\ot b+(-1)^{|a||b|}a\ot [[x, y], b]\eqno(3.6)$$ defines a Leibniz superalgebra structure on the vector space $\mg\ot\mg$
(see [KP] for that in Leibniz algebras case).

3. The {\it general linear Leibniz superalgebra} $\gl(m, n, D)$ is
generated by all $n\times n$ matrices with coefficients from a
dialgebra $D$, and $m, n\ge 0, n+m\ge2$ with the bracket
$$[E_{ij}(a), E_{kl}(b)]=\de_{jk}E_{il}(a\vdash b)-(-1)^{\tau_{ij}\tau_{kl}}\de_{il}E_{kj}(b\dashv a),\eqno(3.7)$$ for all $a, b\in D$.

Clearly, $\gl(m, n, D)$ is a Leibniz superalgebra. If $D$ is an associative superalgebra, then  $\gl(m, n, D)$ becomes a Lie superalgebra.

By definition, the {\it special linear Leibniz superalgebra} with
coefficients in $D$ is
$$\fsl(m, n, D):=[\,\gl(m, n, D), \gl(m, n, D)\,].$$ Notice that if $n\ne m$ the Leibniz superalgebra $\fsl(m, n, D)$ is simple.

The special linear Leibniz superalgebra $\fsl(m, n, D)$ has generators
$E_{ij}(a), 1\le i\ne j\le m+n, a\in D$, which satisfy the
following relations:
\begin{eqnarray*}
&&[E_{ij}(a), E_{kl}(b)]=0,\  \hbox{if}\ i\ne l,\ \hbox{and}\  j\ne k;\\
&&[E_{ij}(a), E_{kl}(b)]=E_{il}(a\vdash b),\  \hbox{if}\ i\ne l,\ \hbox{and}\  j= k;\\
&&[E_{ij}(a), E_{kl}(b)]=-(-1)^{\tau_{ij}\tau_{kl}}E_{kj}(b\dashv a),\
\hbox{if}\ i= l,\ \hbox{and}\  j\ne k,
\end{eqnarray*}

4.  The {\it Steinberg Leibniz superalgebra} $\stl(m, n, D)$
(\cite{L}) is a Leibniz superalgebra generated by symbols
$u_{ij}(a)$, $1\le i\ne j\le n$, $a\in D$, subject to the
relations
\begin{eqnarray*}
&&v_{ij}(k_1a+k_2b)=k_1v_{ij}(a)+k_2v_{ij}(b),\  \hbox{ for } \ a, b\in D, \ k_1, k_2\in K;\\
&&[v_{ij}(a), v_{kl}(b)]=0,\  \hbox{ if }\ i\ne l,\ \hbox{ and } \  j\ne k;\\
&&[v_{ij}(a), v_{kl}(b)]=v_{il}(a\vdash b),\  \hbox{ if } \ i\ne l,\ \hbox{and } \  j= k;\\
&&[v_{ij}(a), v_{kl}(b)]=-(-1)^{\tau_{ij}\tau_{kl}}v_{kj}(b\dashv a),\  \hbox{ if } \ i= l,\ \hbox{and} \  j\ne k,
\end{eqnarray*}
where $1\le i\ne j\le m+n$, $a\in D$. It is clear that the last
two relations make sense only if $m+n\ge3$. See \cite{L} and
\cite{HL} for more details about the Steinberg Leibniz
superalgebra.

We also denote by {\bf SLeib} and {\bf  SLie} the categories of
Leibniz superalgebras and Lie superalgebras over $K$ respectively.

For any super dialgebra $D$, define
$$[x, y]=x\vdash y-(-1)^{|x||y|}y\dashv x,\eqno(3.8)$$
then $D$ equipped with this
bracket is a Leibniz superalgebra. We denoted it by $ D_L$.
The canonical morphism $D\to D_{L}$ induces a functor
$(-):$ {\bf SDias}$\to${\bf SLeib}.

\noindent{\bf Remark.} For a super dialgebra $D$, if we define
$$[x, y]=x\dashv y-(-1)^{|x||y|}y\vdash x,\eqno(3.9)$$ then $(D,\,
[\, , ])$ is a right Leibniz superalgebra.

For a Leibniz superalgebra $L$, let $L_{LS}$ be the quotient of
$L$ by the ideal generated by elements $[x, y]+(-1)^{|x||y|}[y,
x]$, for all $x, y\in L$. It is clear that $ L_{LS}$ is a Lie
superalgebra. The canonical epimorphism $: L\to L_{LS}$ is
universal among the maps from $L$ to Lie superalgebras. In other
words the functor $(-)_{LS}:$ {\bf SLeib}$\to${\bf SLie} is left
adjoint to $inc:$ {\bf SLie}$\to${\bf SLeib}.

Moreover we have the following commutative diagram of categories and functors
\begin{eqnarray*}
{\bf SDias}&\stackrel{-}{\to}& {\bf SLeib}\\
\downarrow&&\downarrow\\
{\bf SAss}&\stackrel{-}{\to}& {\bf SLie}
\end{eqnarray*}

As in the Leibniz algebra case, the {\it universal enveloping
super dialgebra} of a Leibniz superalgebra $L$ is
$$Ud(L):=(T(L)\ot L\ot T(L))/\{[x, y]-x\vdash y+(-1)^{|x||y|}y\dashv x|x, y\in L\}.$$

\begin{prop}
 The functor $Ud:{\bf SLeib}\to {\bf SDias}$ is left adjoint to the functor $-:{\bf SDias}\to {\bf SLeib}$ (see \cite{LP}). \hfill $\rule[-.23ex]{1.0ex}{2.0ex}$
\end{prop}

Let $V$ be a $\z_2$-graded $K$-vector space. The {\it free Leibniz
superalgebra} $\frak L(V)$ is the universal Leibniz superalgebra
for maps from $V$ to Leibniz superalgebras.

Let $\bar T(V):=\oplus_{n\ge1}V^{\ot n}$ be the reduced tensor module. Then we have

\begin{lemm}
 $\bar T(V)$ is the free Leibniz superalgebra over $V$ with the bracket defined inductively by

$(1)$ $[v, x]=v\ot x, \quad \forall x\in \bar T(V), \; v\in V,$

$(2)$ $[y\ot v, x]=[y, v\ot x]-(-1)^{|y||v|}v\ot [y, x]$,  if $x, y\in \bar T(V)$ and $v\in V.$
\end{lemm}
\noindent{\bf Proof.} By direct calculation we can check that
$\bar T(V)$ is a Leibniz superalgebra.

For any Leibniz superalgebra $\mg$ and $\phi: V\to \mg$, define
$f:\bar T(V)\to \mg$ inductively by
$$f(v)=\phi(v)\quad \hbox{and}\quad f(v_1\ot\cdots\ot v_n)=[f(v_1\ot\cdots\ot v_{n-1}), f(v_n)],$$
where the latter is the bracket in $\mg$. Note this definition is
forced by relation (1) . Since $\mg$ is a Leibniz superalgebra,
$f$ satisfies (2). \hfill $\rule[-.23ex]{1.0ex}{2.0ex}$

Let $L$ be a Leibniz superalgebra. We call a $\z_2$-graded space $M=M_{\bar0}\oplus M_{\bar1}$ a module over $L$ if there are two bilinear maps:
$$[-,-]: L\times M\to M \quad \hbox{and} \quad [-,-]: M\times L\to M$$
satisfying the following three axioms
$$[[x, y], m]=[x, [y, m]]-(-1)^{|x||y|}[y, [x, m]],\leqno(\hbox{SLLM})$$
$$[[x, m], y]=[x, [m, y]]-(-1)^{|x||m|}[m, [x, y]],\leqno(\hbox{SLML})$$
$$[[m, x], y]=[m, [x, y]]-(-1)^{|m||x|}[x, [m, y]],\leqno(\hbox{SMLL})$$
for any $m\in M$ and $x,\, y\in L$.

Given a Leibniz superalgebra $L$, let $C^n(L, M)$ be the space of
all super skew-symmetric $K$-linear homogeneous mapping
$L^{\otimes n}\to M,\ n\ge0$ and $C^0(L, M)=M$. Let
$$d^n:C^n(L, M)\to C^{n+1}(L, M)$$
be an $F$-homomorphism defined by
\begin{eqnarray*}
&&(d^nf)(x_1, \dots, x_{n+1})\\
&:=&\sum_{i=1}^{n}(-1)^{(|f|+|x_1|+\cdots+|x_{i-1}|)|x_i|}(-1)^{i-1}[x_i, f(x_1, \cdots, \hat x_i, \cdots, x_{n+1})]\\
&-&(-1)^{(|f|+|x_1|+\cdots+|x_n|)|x_{n+1}|}(-1)^{n}[f(x_1, \cdots, x_{n}), x_{n+1}]\\
&+&\sum_{1\le i<j\le n+1}(-1)^{|x_i|(|x_{i+1}|+\cdots+|x_{j-1}|)}(-1)^{i}f(x_1, \cdots, \hat x_i, \cdots, x_{j-1}, [x_i, x_j], x_{j+1}, \cdots, x_{n+1}),
\end{eqnarray*}
where $f\in C^n(L, M)$.

From \cite{D} we see that $$d^{n+1}d^n=0, \quad  n\ge0.$$
Therefore, $(C^*(L, M), d)$ is a cochain complex, whose cohomology
is called the cohomology of the Leibniz algebra $L$ with
coefficients in the representation $M$:
$$HL^*(L, M)=H^*(C^*(L, M), d).$$

Suppose that $L$ is a Leibniz superalgebra over $K$. For any $z\in
L$, we define $\ad z\in \hbox{End}_kL$ by
$$\ad z(x)=[z, x], \quad\forall x\in L.\eqno(3.10)$$
It follows (2.1) that
$$\ad z([x, y])=[\ad z(x), y]+(-1)^{|z||x|}[x, \ad z(y)]\eqno(3.11)$$
for all $x, y\in L$. This says that $\ad z$ is a super derivation of degree $|z|$
of $L$. We also call it the inner derivation of $L$.

\noindent{\bf Remark.} In right Leibniz superalgebras case, (3.11)
also holds if we define $$\ad z(x)=-(-1)^{|x||z|}[x, z],
\quad\forall x\in L.\eqno(3.12)$$

Similarly we also have the definition of general super derivation
of a Leibniz superalgebra $L$. By definition a super derivation
of degree $s, s\in \z_2$ of $L$ is an endomorphism
$\mu\in\hbox{End}_sL$ with the property
$$\mu([a, b])=[\mu(a), b]+(-1)^{s|a|}[a, \mu(b)].$$
We denote by Inn($L$), Der$(L)$ the sets of all inner derivations,
super derivations of L respectively. They are also Leibniz
superalgebras.

For a Lie superalgebra $L$, $HL^1(L, M)=H^1(L, M)=\der(L,
M)/\hbox{Inn}(L, M)$.

\section{ UNIVERSAL CENTRAL EXTENSIONS OF LEIBNIZ SUPERALGEBRAS}

As in the cases of Lie algebras and Lie superalgebras and Leibniz
algebras, the universal central extensions of Leibniz super
algebras also play an important role in the theory of Leibniz
superalgebras.

\subsection{Central extensions}

\begin{defi}
A central extension of a Leibniz superalgebra $L$ is a short exact
sequence in the category {\bf SLeib}:
$$0\to Z\to \hat L\stackrel{\a}{\to} L\to 0,\eqno(4.1)$$
where $Z$ is in the center of $\hat L$ and $Z$ is called the
kernel of the central extension (4.1).
\end{defi}
We also call a Leibniz superalgebra $L$ {\it perfect} if $[L,
L]=L$ and sometimes denote the above central extension by $\a:\hat
L\to L$.

Two central extensions
$$0\to C\to \hat L\to L\to 0$$
$$0\to C\to \hat L'\to L\to 0$$ of $L$ are said to be {\it equivalent} if there exists an isomorphism of Leibniz superalgebras from
$\hat L\to \hat L'$ such that the following diagram
\begin{eqnarray*}
0\to C\to &\hat L&\to \,L\to 0\\
||\quad &\downarrow j&\ \quad ||\\
0\to C\to &\hat L'&\to L\to 0
\end{eqnarray*}
commutes.

The set of equivalence classes of such central extensions is known
to be parameterized by the second cohomology group $HL^2(L, C)$.
To be precise, we first introduce $Z^2(L, C)$ and $B^2(L, C)$ as
follows. Set
$$Z^2(L, C)=\{\psi:L\times L\to C\mid \psi([a, b], c)=\psi(a, [b, c])-(-1)^{|a||b|}\psi(b, [a, c]), \quad \forall\; a,\, b,\, c\in L \}\eqno(4.2)$$ and
$$B^2(L, C)=\{f:L\times L\to C\mid f(a, b)=g([a, b])\ \hbox{for some}\ K\hbox{-linear map}\ g:L\to C  \}.\eqno(4.3)$$
An element in $Z^2(L, C)$ is called a Leibniz super 2-cocycle.

\begin{lemm}
Let $C$ be a $K$-module. The second cohomology group $HL^2(L,
C)=Z^2(L, C)/B^2(L, C)$ is in one-to-one correspondence with the
set of the equivalence classes of central extensions of $L$ by
$C$.
\end{lemm}

\noindent{\bf Proof.} The proof is classical and essentially the
same as in the case of Lie superalgebras. \hfill
$\rule[-.23ex]{1.0ex}{2.0ex}$

\subsection{Universal central extensions}

Now we define the universal central extension of a Leibniz superalgebra $L$.

\begin{defi}
The central extension (4.1) of $L$ is called the universal central
extension (UCE) if the following conditions hold.

1. $\hat L$ is perfect.

2. For any central extension $\b:\hat L'\to L$, there exists $\r: \hat L\to\hat L'$ such that $\a=\b\cdot\r$.
\end{defi}

\noindent{\bf Remark.} A UCE is unique up to isomorphism of
Leibniz superalgebras if it exists.

We notice that the mapping $\r$ in the Definition 4.3 is unique.
In fact we have:

\begin{lemm}
Let $(X, f)$ and $(Y, g)$ be central extensions of a Leibniz
superalgebra $L$. If $Y$ is perfect, then there exists only one
homomorphism $h$ from $Y$ to $X$ such that $f\circ h=g$.
\end{lemm}
\noindent{\bf Proof.} Suppose that there are two homomorphisms
$h_i:Y\to X$ such that $f\circ h_i=g, i=1, 2$.

For $a, b\in Y$, let $h_i(a)=a_i\in X, h_i(b)=b_i\in X, i=1, 2$. Then $f(a_1)=f(h_1(a))=g(a)=f(h_2(a))=f(a_2)$. Similarly, $f(b_1)=f(b_2)$. Therefore $a_1=a_2+r_a, b_1=b_2+r_b$, where $r_a, r_b\in\ker f$. Hence
$$h_1([a, b])=[h_1(a), h_1(b)]=[a_1, b_1]=[a_2+r_a, b_2+r_b]=[a_2, b_2]=[h_2(a), h_2(b)]=h_2([a, b]).$$ Since $Y$ is perfect, $h_1=h_2$.  \hfill $\rule[-.23ex]{1.0ex}{2.0ex}$

\begin{prop}
A Leibniz superalgebra $L$ admits a universal central extension
$\hat L$ if and only if $L$ is perfect.
\end{prop}
\noindent{\bf Proof.} Suppose that $\a: \hat L\to L$ is the UCE.
By definition, $\hat L$ is perfect, and hence
$$L=\a(\hat L)=\a([\hat L, \hat L])=[\a(\hat L), \a(\hat L)]=[L, L].$$

Next let us suppose that $L$ is perfect. We set
$$W=(L\ot L)/I,$$ where $I$ is the ideal generated by $[a, b]\ot c-a\ot[b, c]+(-1)^{|a||b|}b\ot[a, c]$ for all $a, b, c\in L$. Let $w: L\ot L\to W$ be the canoical projection. It is clear that $w\in Z^2(L, W)$. We consider the central extension
$$0\to W\to L_w\to L\to 0$$ defined by $w$. Using this central extension, we construct the universal central extension of $L$. Let $V$ be an arbitrary $K$-module and $f\in Z^2(L, V)$. Since
$$f([x, y], z)-f(x, [y, z])+(-1)^{|x||y|}f(y, [x, z])=0,\eqno(4.4)$$
we have a $K$-linear map $\psi':W\to V$ such that $w(x, y)\to f(x,
y)$.

Let us define $\phi':L_w\to L_f$ by
$$\phi'((x, u))=(x, \psi'(u)).$$
It is clear that  $$\a=\b\cdot\phi'.$$
Now we set $\hat L=[L_w, L_w]$. Since $L$ is perfect, it follows that $\hat L+W=L_w$. This implies that $\hat L$ is perfect since
$$\hat L=[\hat L+W, \hat L+W]=[\hat L, \hat L].$$
Furthermore if we set
$$C=W\cap \hat L,$$ then we have a central extension
$$0\to C\to \hat L\to L\to 0$$
such that $\hat L$ is perfect.

Now if we define $\phi$ as the restriction of $\phi'$ to the
subalgebra $\hat L$, then we have $$\a|_{\hat L}=\b\cdot\phi.$$
Therefore, $\hat L\to L$ is the universal central extension of
$L$. \hfill $\rule[-.23ex]{1.0ex}{2.0ex}$

The following two propositions are clear (see \cite{BM} and
\cite{P} for the Lie algebra case, see \cite{LL} for the Leibniz
algebra case).

\begin{prop}
If $L$ is a perfect Leibniz superalgebra and $\hat L$ is a
universal central extension of $L$, then every derivation of $L$
lifts to a derivation of $\hat L$. If $L$ is centerless, the lift
is unique and $\der(\hat L)= \der(L)$.
\end{prop}

\noindent{\bf Proof.} We use the notation as the proof Proposition
4.5 and may assume that $\hat L=[L\oplus W,  L\oplus W]$. Let
$\mu\in \der(L)$. Then $\mu$ induces a linear mapping on $L\ot L$
in the usual way:
$$\mu(x\ot y)=\mu(x)\ot y+(-1)^{|x||\mu|}x\ot \mu(y).$$ This mapping stabilizes $I$, so it induces a map, also
denoted
by $\mu$, on $W$.

Thus $x+u\mapsto \mu(x)+\mu(u)$ defines a derivation on $L_w$
whose restriction to $L$ is the required lifting.

Suppose that the kernel of the central extension is $C$. Since $L$
is centerless, the center of $\hat L$ is $C$, and every derivation
of $L$ induces a derivation of $L\cong \hat L/C$. If $\mu_1$ and
$\mu_2$ are derivations of $\hat L$, both lifts of $\mu$, then for all
$x\in \hat L$, $\mu_2(x)=\mu_1(x)+a(x)$ where $a(x)\in C$. Then
calculation

$$\mu_2([x, y])=[\mu_2(x), y]+[x, \mu_2(y)]=[\mu_1(x), y]+[x, \mu_1(y)]=\mu_1([x, y])$$ shows that $\mu_2=\mu_1$. It follows that $$\der(\hat L)=\der L.\quad  \rule[-.23ex]{1.0ex}{2.0ex}$$

\begin{prop}
Let $L$ be a perfect Leibniz algebra and $\hat L$ be its universal
central extension. Every automorphism $\theta$ of $L$ admits a
unique extension $\hat\theta$ to an automorphism $\hat L$.
Furthermore, the map $\theta\to\hat\theta$ is a group
monomorphism. If $L$ is centerless, then $\aut L\cong\aut\hat L$.
\end{prop}

\noindent{\bf Proof.} We use again the notation as in the proof of
Proposition 4.5 and may assume that $\hat L=[L\oplus W,  L\oplus
W]$. Every automorphism $\theta$ induces an automorphism
$\theta_W$ of $W$ via $\theta_W(x\ot y)=\theta(x)\ot\theta(y)$. It
is clear that $\theta$ extends to an automorphism $\theta_{L_w}$
of $L_w$ satisfying $\theta_{L_w}: x+u\to
\theta(x)+\theta_{W}(u)$. By restriction, $\theta_{L_w}$ induces
an automorphism $\hat\theta$ of $\hat L$.

It is clear from the definition that $\xi:\theta\mapsto \hat\theta$ is a group homomorphism.
Suppose $\hat\theta=1$. Then for all $x, y\in L\subset L_w$ we have
\begin{eqnarray*}
[x, y]+\psi(x, y)&=&[x, y]_{L_w}=\hat\theta([x, y]_{L_w})=\theta_{L_w}([x, y]_{L_w})=[\theta_{L_w}(x), \theta_{L_w}(y)]_{L_w}=[\theta(x), \theta(y)]_{L_w}\\
                 &=&[\theta(x), \theta(y)]+\psi(\theta(x), \theta(y))=\theta([x, y])+ \psi(\theta(x), \theta(y)).
\end{eqnarray*}
Thus $\theta$ is the identity on $[L, L]=L$. It follows that $\xi$
is injective. Next we show that the lifting of $\theta$ to $\hat
L$ is unique ( and hence equal to $\hat\theta$). Let $\theta_1$
and $\theta_2$ be two lifts of $\theta$ to $\hat L$. Then for
$x\in \hat L$ we have $\theta_1(x)=\theta_2(x)+x_C$ for some
$x_C\in C$, where $C$ is the kernel of the universal central
extension. Thus for all $x, y\in \hat L$ we have
$$\theta_1([x, y]_{\hat L})=[\theta_2(x)+x_C, \theta_2(y)+y_C]_{\hat L}=[\theta_2(x), \theta_2(y)]_{\hat L}=\theta_2([x, y]_{\hat L}).$$
Thus $\theta_1=\theta_2$ on $[\hat L, \hat L]=\hat L$.

For a subspace $E$ of $C$, we set $$\aut(L,
E):=\{\theta\in\aut(L)\mid \theta_W(E)=E\}.\eqno(4.5)$$ If
$\theta\in\aut(L, E)$, then $\hat\theta$ induces an automorphism
of $\hat L/E$ that we will denote by $\hat\theta_E$. With similar
arguments we can show that $\theta\to\hat\theta_E$ is a group
monomorphism from $\aut(L, E)$ into $\aut(\hat L)/E$. Now we show
that this map is surjective if $L$ is centreless.

Let $\sg\in \aut(\hat L)/E$, we know that $\sg$ admits a unique
extension $\hat\sg$ to $\hat L$ since $\hat L$ is the UCE of $\hat
L/E$. If $L$ is centreless, the center of $\hat L/E$ is $C/E$. In
this case $\sg\in\aut(\hat L)/E$ induces an automorphisms of
$(\hat L/E)/(C/E)\cong L$ that will be denoted by $\theta$. Then
both $\hat\theta$ and $\hat\sg$ are extensions of $\theta$ to
$\hat L$, and therefore equal. Now $\hat\sg$ stabilizes $E$(since
$E$ is the kernel of $\hat L\to \hat L/E$) and therefore so does
$\hat\theta$. But the restriction of $\hat\theta$ to $C\subset W$
coincides with $\theta_W$. This show that $\theta$ is in fact an
element of $\aut(L, E)$. By taking $E=\{0\}$, we obtain
$$\aut(\hat L)=\aut(L).\eqno \rule[-.23ex]{1.0ex}{2.0ex}$$

\section{UNIVERSAL CENTRAL EXTENSIONS OF THE BASIC CLASSICAL LIE SUPERALGEBRAS}

Throughout this section $\mg$ always denotes a basic classical
simple Lie superalgebra and $D$ denotes a unital commutative
dialgebra over $\c$ and $A=\c[t_1^\pm,\cdots,t_\nu^\pm](\nu\ge1)$.
For convenience, we always assume that $\mg$ is not of type $A(n,
n)$.

In this section we shall determine the universal central extension
of the Lie superalgebra ${\mg}\ot D$. Now we describe the main
theorem of this section.

\begin{theo}
The universal central extesion of ${\mg}\ot D$ is
$$\wg=\mg\ot D\oplus\Omega_D^1 $$ with the bracket
$$[X\ot a, Y\ot b]=[X, Y]\ot (a\vdash b)+(X, Y)b\dashv da, \quad, X, Y\in\wg, a, b\in D,$$ where $( , )$ is an even supersymmtric invariant bilinear form on $\mg$ describe in section 3.
\end{theo}

In what follows, we mainly prove Theorem 5.1.

First we introduce some notation. Let $$0\to Z\to \frak u\stackrel{\pi}{\to} \frak a\to 0$$ be a central extension. Notice that, for $x, y\in \frak a$ and $x', y'\in \frak u$ such that $\pi(x')=x, \pi(y')=y$, the commutator $[x', y']$ does not depend on the choice of the inverse images $x'$ and $y'$. Hence we denote $[x', y']$ by $[x, y]'$

\begin{prop}
Suppose $\frak g$ is a basic classical Lie superalgebra which is
not of type $A(n, n)$. Let $Z$ be a free $\c$-module and
$$0\to Z\to \mg'(D)\stackrel{\pi}{\to} \mg\ot D\to 0$$ a central extension of $\mg\ot D$. Then the bracket of $\mg'(D)$ can be described as
$$[(X\ot a)', (Y\ot b)']=([X, Y]\ot(a\vdash b))'+(X, Y)\{a, b\},\eqno(5.1)$$
where $( , )$ is a non-degenerate even supersymmetric invariant bilinear form on $\mg$, $\{ ,\}:D\times D\to Z$ satisfies
$$\{a\vdash b, c\}=\{a, b\vdash c\}+\{b, a\vdash c\}.\eqno(5.2)$$
\end{prop}

\noindent{\bf Proof.}  (1) We first consider the case of $\fsl_2$,
which is regard as a trivial Lie superalgebra, so its proof is
same as that in Leibniz central extension case given in \cite{LL}
and [Gao1]. It is also similar to Lemma 4.12 in \cite{IK} (just
delete $\{a, b\}=-\{b, a\}$ and replace $ab$ by $a\vdash b$).

(2) Now we consider the case of $\frak{osp}(1, 2)$:
$$\frak{osp}(1, 2)=\c X_+\oplus\c x_+\oplus\c H\oplus\c x_-\oplus\c X_-,$$ with
$$[H, X_{\pm}]=\pm4X_{\pm}, \  [X_+, X_-]={1\over 2}H,\  X_\pm=\pm{1\over4}[x_\pm, x_\pm], $$
$$[H, x_{\pm}]=\pm2x_{\pm}, \  [x_+, x_-]=H,\  [X_\pm, x_\pm]=-x_\pm. $$

With this Chevalley basis\ ($|x_\pm|=\bar1$, $|X_\pm|=|H|=\bar0$),
we can give a non-dengenerate even supersymmetric bilinear form on
$\frak{osp}(1, 2)$ as follows:
$$(H, H)=2,\quad (x_+, x_-)=1, \quad (X_+, X_-)={1\over4}.$$

We set
\begin{eqnarray*}
(X_{\pm}\ot a)':&=&\pm{1\over 4}[H\ot 1, X_{\pm}\ot a]\\
(x_{\pm}\ot a)':&=&\pm{1\over 2}[H\ot 1, x_{\pm}\ot a]\\
(H\ot a)':&=&[x_+\ot 1, x_-\ot a]'\\
\{a, b\}:&=& {1\over2}[(H\ot a)', (H\ot b)'].
\end{eqnarray*}
We shall show these elements satisfy (5.2). It is suffice to show
the following formulas:

F1) $[(H\ot a)', (Y\ot b)']=\b(H)(Y\ot(a\vdash b))'$ and $[(Y\ot a)', (H\ot b)']=-\b(H)(Y\ot(a\vdash b))'$, where $Y\in \mg_\b$.

F2) $[(x_{\pm}\ot a)', (x_\mp\ot b)']=\pm(H\ot(a\vdash b))'+\{a, b\}$ and $[(H\ot a)',\ (H\ot b)']=2\{a, b\}$.

F3) $[(x_\pm\ot a)', (x_\pm\ot b)']=\pm4(X_\pm\ot(a\vdash b))'$.

F4) $[(X_\pm\ot a)', (x_\mp\ot b)']=-(x_\pm\ot(a\vdash b))'$ and $[(x_\pm\ot a)', (X_\mp\ot b)']=(x_\pm\ot(a\vdash b))'$.

F5) $[(X_\pm\ot a)', (X_\mp\ot b)']=\pm{1\over2}(H\ot(a\vdash b))'+{1\over4}\{a, b\}$.

F6) $[(X_\pm\ot a)', (x_\pm\ot b)']=0=[(x_\pm\ot a)', (X_\pm\ot b)'], \   [(X_\pm\ot a)', (X_\pm\ot b)']=0$.

We prove F1). By definition and Jacobi identity
\begin{eqnarray*}
& &[(H\ot a)', (X_\pm\ot b)']\\
&=&[(H\ot a)', \pm{1\over4}[(H\ot1)', (X_\pm\ot b)']\\
&=&\pm{1\over4}([[(H\ot a)', (H\ot1)'], (X_\pm\ot b)']+[(H\ot 1)', [(H\ot a)', (X_\pm\ot b)']]).
\end{eqnarray*}
Since
\begin{eqnarray*}
&&[(H\ot a)', (H\ot1)']\equiv 0(\mod Z)\\
&&[(H\ot a)', (X_\pm\ot b)']= [H\ot a, X_{\pm}\ot b]'\equiv \pm 4(X_{\pm}\ot(a\vdash b))' (\mod Z).
\end{eqnarray*}

So for $Y=X_{\pm} $, F1) is true. Similar we can prove F1) for $Y=x_{\pm}$.

The formulas F2)---F6) are just follows from direct calculation
similar that in Lemma 4.13 in \cite{IK}.

(3) Suppose $\frak g$ is a basic classical Lie superalgebra which
is not of type $A(n, n)$.

Let $\{X_\a, H_\a|\a\in \Delta\}$ be a Chevalley basis as in
Section 3.

For $\a\in \Delta_{\bar0}$ such that $\fr{1}{2}\a\notin\Delta$, we set
$$(X_{\b}\ot a)':={1\over \b(H_\a)}[H_\a\ot 1, X_{\b}\ot a]'$$
$$(H_{\b}\ot a)':=\sg(\b)[X_\b\ot 1, X_{-\b}\ot a]',$$ for $\b=\pm\a$ and
$$\tilde\mg(\a):=\{\c X_\a\oplus\c H_\a\oplus\c X_{-\a}\}\ot D,$$
$$\tilde\mg'(\a):=\pi^{-1}(\tilde\mg(\a)).$$

For $\a\in\Delta_{\bar1}$ such that $2\a\in\Delta$, we set
$$(X_{\b}\ot a)':={1\over \b(H_\a)}[H_\a\ot 1, X_{\b}\ot a]'$$
$$(H_{\b}\ot a)':=\sg(\b)[X_\b\ot 1, X_{-\b}\ot a]',$$ for $\b=\pm\a, \pm2\a$ and
$$\tilde\mg(\a):=\{\c X_{2\a}\oplus\c X_\a\oplus\c H_\a\oplus\c X_{-2\a}\oplus\c X_{-\a}\}\ot D,$$
$$\tilde\mg'(\a):=\pi^{-1}(\tilde\mg(\a)).$$

Then from Case (2) and Case (3), we see that for $\a\in\Delta$
such that $(\a, \a)\ne0$, there exists $\{ , \}:D\times D\to Z$
such that
$$[(X\ot a)', (Y\ot b)']=([X, Y]\ot(a\vdash b))'+(X, Y)\{a, b\}$$ and
and $$\{a\vdash b, c\}=\{a, b\vdash c\}+\{b, a\vdash c\}.$$

For $\a\in\Delta_{\bar1}$ such that $2\a\notin\Delta$ i.e. $(\a, \a)=0$, we fix $\r\in\Delta$ as in Lemma 3.5, and set
$$(X_{\b}\ot a)':={1\over \b(H_\r)}[H_\r\ot 1, X_{\b}\ot a]'$$
$$(H_{\b}\ot a)':=\sg(\b)[X_\b\ot 1, X_{-\b}\ot a]',$$ for $\b=\pm\a$.

For above $\r$, we also introduce the following notation:
$$\bar\a:=\cases{\a, \quad (\a, \a)\ne0\cr \r, \quad (\a, \a)=0}.$$

Now it is sufficient to prove the following formulas, whose proofs
are similar to Lemma 4.16---4.18, 4.20---4.22 in \cite{IK}.

1) For $\a, \b\in\Delta$, we have
$$[(H_\a\ot a)', (X_\b\ot b)']=\b(H_\a)(X_\b\ot(a\vdash b))'.$$

2) For each $\a\in\Delta$, we set $\{a, b\}_\a=[(H_{\bar\a}\ot a)',\ (H_\a\ot b)']$ and $\{a, b\}=\fr{\sg_\a}{(X_\a, X_{-\a})}\{a, b\}_\a$, where $\sg_\a$ is defined as before. Then $\{ ,\}$ is independent of the choice of $\a$ and
$$[(X_\a\ot a)', (X_{-\a}\ot b)']=\sg_\a(H_\a\ot(a\vdash b))'+(X_\a, X_{-\a})\{a, b\},$$
$$[(H_\a\ot a)', (H_{\b}\ot b)']=(H_\a, H_\b)\{a, b\}.$$

3) For $\a, \b\in\Delta$ such that $\a+\b\ne0$.
$$[(x_\a\ot a)', (x_\b\ot b)']=N_{\a,\b}(X_{\a+\b}\ot(a\vdash b))', $$ where $N_{\a,\b}$ is defined as in section 3.

4) For $\a, \b\in\Delta$ such that $\a+\b\in\Delta$, we have
$$(H_{\a+\b}\ot a)'={\vep_{\a+\b}\over\vep_\a}(H_\a\ot a)'+{\vep_{\a+\b}\over\vep_\b}(H_\b\ot a)'.$$

5) $\{ ,\}:D\times D\to Z$ satisfies
$$\{a\vdash b, c\}=\{a, b\vdash c\}+\{b, a\vdash c\}.\eqno \rule[-.23ex]{1.0ex}{2.0ex}$$

Let us complete the proof of Theorem 5.1 by using Proposition 5.2.
Let $\phi: D\times D\to \Omega_D^1$ defined by $(a, b)\to b\dashv
da$. Since the Leibniz superalgebra $\mg\ot D$ is perfect, it has
a universal central extension by Proposition 4.1. By definition of
$\Omega_D^1$ and Proposition 5.2, for any central extension
$$0\to Z\to \mg'(D)\to \mg\ot D\to 0$$
there exists $\psi:\Omega_D^1\to Z$ such that
$$\{ , \}=\psi\cdot\phi.$$
Hence the central extension
$$0\to \Omega_D^1\to \tilde\mg\to \mg\ot D\to 0$$ is the universal central extension.\hfill $\rule[-.23ex]{1.0ex}{2.0ex}$

\begin{coro}
Under the assumption of Theorem 5.1, we have
$$HL_2(\mg\ot D)\cong \Omega_D^1.$$
\end{coro}

Especially, if $D$ degenerate a unital cummutative algebra
$A=\c[t_1^\pm,\cdots,t_\nu^\pm](\nu\ge1)$, we have
\begin{coro}
The UCE of $\mg\ot A$ in the category of Leibniz superalgebras
is
$$\wg=\mg\ot A\oplus\Omega_A^1 $$ with the bracket
$$[X\ot a, Y\ot b]=[X, Y]\ot (ab)+(X, Y)bda, \quad X, Y\in\wg,\  a, b\in A,$$ where $( , )$ is an even supersymmtric invariant bilinear form on $\mg$ described in section 3.
\end{coro}

With Theorem 4.7 in \cite{IK}, the methods of (4.6) in [LP] and
Corollary 5.5 we also have

\begin{coro}
Let $\mg(A)$  be the UCE of $\mg\ot A$ in the category of Lie
algebras (see section 1) and $\tilde\mg$ as in Corollary 5.4. Then
the UCE of $\mg(A)$ in the category of Leibniz superalgebras is
$\wg$ with kernel $\im B$, where $B$ is the modified Connes
operator defined in \cite{HL} (also see [Lo3]).
\end{coro}

With Corollary 5.4 and Proposition 4.5 and 4.6 we have
\begin{coro}
The derivation algebra of the Leibniz superalgebra $\tilde\mg$ is
$$\der(\tilde\mg)=\der(A)\ltimes{Inn}(\mg)$$ and
the automorphism group of $\tilde\mg$ is
$$\aut(\tilde\mg)=\aut(\mg\ot D),$$ where $\der(A)=\sum_{i=1}^{\nu}A{\pt\over \pt t_i}$ is the derivation algebra of $A$ as an assicoative algebra and $\aut(\mg\ot A)$ is the automorphism group of Lie superalgebras $\mg\ot A$. \hfill $\rule[-.23ex]{1.0ex}{2.0ex}$
\end{coro}

\noindent{\bf Remarks.} 1. If we consider $\mg$ is the type $A(n,
n)$, we also show that $\wg:=\bar\mg\ot D\oplus\Omega_D^1$, where
$\bar\mg$ defined as in section 1, is also the universal central
extension of $\mg$, i.e.
\begin{prop}
If $\mg$ is a basic classical Lie superalgebra and $D$ is a unital commutative dialgebra, then we have
$$HL_2(\mg\ot D)\cong \cases {\Omega_D^1,\ \hbox{if}\ \mg\  \hbox{is not of type}\ A(n, n),\cr \Omega_D^1\oplus D^{\oplus3}, \ \hbox{if}\ \mg\  \hbox{is of type}\ A(1, 1),\cr
\Omega_D^1\oplus D\ \hbox{if}\ \mg\  \hbox{is of type}\ A(n, n),
n>1.}$$ (see \cite{IK}).
\end{prop}
2. For a general associative algebra $A$, the universal central
extension of the Lie superalgebra $\fsl(m, n)\ot A$ in the Leibniz
superalgebras category is obtained in \cite{HL}.

\begin{theo}$\cite{HL}$
Let $k$ be a field\,(char $k\ne 2, 3$) and $A$ an associative and
unital $k$-algebra. For $m+n\ge 3,\, m, n\ge0$, the universal
extension of $\fsl(m, n, A):=\fsl(m, n)\ot A$ in the category of
Leibniz superalgebras is
$$0\to HH_1(A)\to \frak{stl}(m, n, A)\to \fsl(m , n)\ot A\to 0.$$
In particular, if $A$ is commutative, then
$$HL^2(\fsl(m, n)\ot A)\cong\Omega_{A|k}^1,$$ where $\frak{stl}(m, n, A)$ is the special Steinberg Leibniz superalgebra $\stl(m, n, D)$ (section 4) for $D$ degenerated by an associative algebra $A$.
\end{theo}

3. More generally, for a unital dialgebra $D$, the universal
central extension of the Leibniz superalgebra $\fsl(m, n, D)$ is
obtained in \cite{L}.

\begin{theo}$\cite{L}$
Let $k$ be a field\,$($char $k\ne 2, 3$$)$ and $D$ a unital
dialgebra. For $m+n\ge 3,\, m, n\ge0$, the universal extension of
$\fsl(m, n, D)$ is
$$0\to HHL_1(D)\to \frak{stl}(m, n, D)\to \fsl(m , n, D)\to 0.$$
In particular, if $D$ is commutative $($i.e. $a\vdash b=b\dashv a$
for all $a, b\in D$$)$, then
$$HL^2(\fsl(m, n, D))\cong\Omega_{D}^1,$$ where $HHL_1(D)$ is first modified Frabetti's Homology
group for dialgebras defined in \cite{L} $($also see \cite{F1} and
\cite{F2}, $HHL_1(D)\cong \Omega_{D}^1$ if $D$ is commutative$)$.
\end{theo}

4. Let $D$ be a unital super dialgebra, then $\stl(m, n, D)$ and
$\mg\ot D$ ($D$ commutative) are also Leibniz superalgebras
according to (4.1) and (4.7) and (4.8). They play key roles in
studying the Leibniz superalgebras graded by finite root systems
as in the Lie algebra, Leibniz algebra and Lie superalgebra case
(see \cite{BM}, \cite{BE}, \cite{BZ} and \cite{LH2}).

\vskip10pt \centerline{\bf ACKNOWLEDGMENTS}

\vskip5pt The part of this work was done during the study of the
first author in the Department of Mathematics at University of
Bielefeld for his Ph. D. degree and in the Department of
Mathematics at Shanghai Jiaotong University for his Post-doctor
study. He would like to expresses his special gratitude to the
`AsiaLink Project' for financial support. He is also deeply
indebted to Prof.s  C.M. Ringel and C.P. Jiang for their kind
hospitality and continuous encouragement and instruction. The
project is supported by the NNSF (Nos. 10431040, 10271047,
19731004), the TRAPOYT and the FUDP from the MOE of China, the
SYVPST from the SSTC, and the Shanghai Priority Academic
Discipline from the SEC, the China Postdoctoral Science
Foundation. Authors give their special thanks to Prof. Y. Gao for
the crucial comments, to Referee for his useful comments and
corrections of many English mistakes.

\vskip10pt
\def\refname{\cen{\normalsize\bf REFERENCES}}

\end{document}